\documentstyle[12pt]{article}
\newlength{\abstractwidth}
\renewcommand{\thefootnote}{\fnsymbol{footnote}}
\renewcommand{\thanks}[1]{\footnote{#1}} 
\newcommand{\starttext}{ \setcounter{footnote}{0}
\renewcommand{\thefootnote}{\arabic{footnote}}}

\newcommand{\be}{\begin{equation}}
\newcommand{\bea}{\begin{eqnarray}}
\newcommand{\eea}{\end{eqnarray}} \newcommand{\ee}{\end{equation}}
 
 \def\ba{\begin{eqnarray}}
\def\ea{\end{eqnarray}}

\def\o{\omega}

\def\log{\,{\rm log}\,}

\def\o{\omega}

\def\o{\omega}

\def\p{\partial}

\def\ddb{{\partial\bar\partial}}

\def\[{{\bf [}}
\def\]{{\bf ]}}

\begin{document}
\starttext \baselineskip=18pt \setcounter{footnote}{0}
\newtheorem{theorem}{Theorem}
\newtheorem{lemma}{Lemma}
\newtheorem{corollary}{Corollary}
\newtheorem{definition}{Definition}
\newtheorem{conjecture}{Conjecture}
\newtheorem{proposition}{Proposition}

\begin{center}
{\Large \bf GEOMETRIC PARTIAL DIFFERENTIAL EQUATIONS FROM UNIFIED STRING THEORIES
\footnote{Contribution to the proceedings of the ICCM 2018 conference in Taipei, Taiwan.
Work supported in part by the National Science Foundation under grant DMS-12-66033 and by the Galileo Galilei Institute of Theoretical Physics.}}

\medskip
\centerline{Duong H. Phong}

\medskip

\begin{abstract}

An informal introduction to some new geometric partial differential equations motivated by string theories is provided. Some of these equations are also interesting from the point of view of non-K\"ahler geometry and the theory of non-linear partial differential equations. In particular, a survey is given of joint works of the author with Teng Fei, Bin Guo, Sebastien Picard, and Xiangwen Zhang.

\end{abstract}

\end{center}

\baselineskip=15pt
\setcounter{equation}{0}
\setcounter{footnote}{0}

\section{Introduction}
\setcounter{equation}{0}

The laws of nature at its most fundamental level have long been a source of inspiration for the theory of geometric partial differential equations. Each of the equations for the four known forces of nature, namely Maxwell's equations for electromagnetism, the Yang-Mills equations for the weak and the strong interactions, and Einstein's equation for gravity, has led to a deep and rich theory, with wide ramifications in many branches of mathematics, including topology and geometry. In the mid 1980's, string theories burst upon the scene as the only consistent candidates for a unified theory of all forces including gravity \cite{GSW}. In their simplest form, the equations from string theories unexpectedly brought in complex structures \cite{CHSW}, together with a completely new motivation for equations from K\"ahler geometry which had only been studied before as generalizations of the Uniformization Theorem \cite{Y}. This confluence of complex geometry with string physics has contributed to some of the most exciting developments in mathematics and theoretical physics in the last few decades, notably mirror symmetry.

\smallskip
However, for many reasons, including phenomenology, the moduli-space problem, and non-perturbative effects such as branes and duality
\cite{BBS, Ki, Po}, it has proven to be useful to consider solutions of string equations which are more general than the simplest one first considered in \cite{CHSW}. A frequent feature of these more general solutions is the incorporation of non-vanishing fluxes, which correspond on the geometric side to Hermitian connections with non-vanishing torsion. Thus physical considerations lead, perhaps surprisingly, to the realm of non-K\"ahler geometry. This is a wide open area in mathematics, still largely unexplored, and natural questions originating from theoretical physics should provide a very valuable guidance.

\smallskip
The equations required for the simplest solutions of string theories had actually been already solved by Yau \cite{Y} and Donaldson \cite{D} and Uhlenbeck-Yau \cite{UY} at the time they were proposed in \cite{CHSW}. But much less is known even at the present time about the more general equations referred to above. Even ignoring the motivation from physics or geometry, they are quite interesting from the point of view of the theory of non-linear partial differential equations. All this is well appreciated in selected circles, but it may not be as widely known as it should be.

\smallskip
The purpose of this lecture is to provide an informal, and necessarily very incomplete, introduction to these equations. The focus will only be on some joint works of the speaker with Teng Fei, Bin Guo, Sebastien Picard, and Xiangwen Zhang, and only some particular equations from $11$-dimensional supergravity, the heterotic string, and the Type II strings also of interest in geometry will be discussed.

\section{Some Aspects of Unified String Theories}
\setcounter{equation}{0}

Of the four known forces of nature, electromagnetism and weak and strong interactions are described by gauge fields, Abelian in the case of electromagnetism, and non-Abelian in the case of the other two forces. A gauge field is a connection $A$ on a vector bundle over space-time, the field strength is given by the curvature $F=dA+A\wedge A$, and the field equations are the Yang-Mills equation, coupled with the Bianchi identity,
\bea
d_A^\dagger F=0,\qquad d_AF=0.
\eea
Gravity is on the other hand described by general relativity, where the field is a metric $g_{ij}$ and the field equation in vacuum is given by Einstein's equation
\bea
R_{ij}=0
\eea
where $R_{ij}$ is the Ricci curvature of the metric $g_{ij}$. While many examples of solutions of the corresponding field equations have been obtained by exploiting a complex structure, the original equations are formulated for a real space-time, which is usually Lorentz, although a Wick rotation to Euclidian is of interest as well. Gauge theories are renormalizable while general relativity is not, and thus a quantum theory of gravity does not appear accessible by standard field theory methods.

\medskip
The unification of all the forces of nature, including gravity, into a single, consistent quantum theory, is one of the grand dreams of theoretical physics. 
So far, practically all attempts with just field theories have run into serious difficulties, and since the mid 1980's, the only known viable candidate has been {\it string theories}, which are theories of one-dimensional extended objects. There are 5 string theories, known as the Type I string, the Heterotic string with either $SO(32)$ or $E_8\times E_8$ as gauge groups, and the Type IIA and Type II B string. The Type I theory is a theory of open strings, the other theories are theories of closed strings. They all take place in a $10$-dimensional space-time \cite{GSW}.
Since the mid 1990's, it has been realized that they are all related by dualities, and should be viewed as manifestations of a new quantum theory, called $M$ theory, which is approximated at low energies by a classical field theory, namely $11$-dimensional supergravity \cite{HW, W, T, BBS, Ki, Po}.

\smallskip
All these theories incorporate supersymmetry, which is a symmetry pairing bosons (represented by tensor fields) with fermions (represented by spinor fields).

\subsection{Supergravity theories}

It is not possible to give here an adequate description of these theories which are difficult and still far from being understood. Instead, we shall just focus on a few mathematical implications which play an important role in the geometric partial differential equations that we shall describe in the sequel.

\smallskip
First, for our purposes, we do not have to deal with the extended objects themselves. Rather, 
in the low-energy limit, string theories and $M$ Theory reduce to just field theories of point particles in a $10$ or $11$-dimensional space-time, and we shall just consider these. Since string theories automatically incorporate gravity and are supersymmetric, their low energy limits are {\it supergravity theories}, in the sense that they are field theories which incorporate gravity and are supersymmetric. The incorporation of gravity means that the field content always includes a metric $G_{MN}$, where $M,N$ are space-time indices.
The requirement of supersymmetry on a higher-dimensional space-time is a severe constraint, and there are very few supergravity theories. In each case, it suffices to specify the bosonic fields, and the fermionic fields and the action are then completely determined by supersymmetry.

\smallskip
We can accordingly describe the field theory limits of string theories and $M$ theory as follows \cite{CJS, BBS, Ki, Po}:

\medskip
$\bullet$ {\it $11$-dimensional supergravity}: 
The bosonic fields are a metric $G=G_{MN}$ and a $4$-form $F=dA_3$ on an $11$-dimensional space-time. The action is 
\bea
\label{11d-action}
I=\int d^{11}x \sqrt{-G}(R-{1\over 2}|F_4|^2)-{1\over 6}\int A_3 \wedge F_4\wedge F_4
\eea
Here $F_4=dA_3$ is the field strength of the potential $A_3$, and all fermionic fields have been set to $0$, as we shall also do systematically below.

\smallskip
$\bullet$ {\it The Type I and Heterotic strings}: 
The bosonic fields are a metric $G_{MN}$, a $2$-form $B_{MN}$, a scalar $\Phi$ (from the gravity multiplet)in $10$-dimensional space-time, and a vector potential $A_M$ (from the vector multiplet),
valued in $SO(32)$ in the case of the Type I string, or in either $SO(32)$ or $E_8\times E_8$ in the case of the Heterotic string. The action is
\bea
\label{Type I}
I=\int d^{10}x \sqrt{-G}(R-|\nabla\Phi|^2-e^{-\Phi}|H|^2-e^{-\Phi/2}Tr[F^2])
\eea
and $F$ is the curvature of $A$, and
\bea
\label{AnomalyCancellation}
H=dB-\omega_{CS}(A)+\omega_{CS}(L),
\eea
where $\omega_{CS}(A)$ is the gauge Chern-Simons form ${\rm Tr}(A\wedge dA-{2\over 3}A\wedge A\wedge A)$, and $\omega_{CS}(L)$ is the Lorentz Chern-Simons form.

\smallskip
$\bullet$ {\it The Type II A and Type II B strings}: the bosonic fields are again the gravity multiplet $(G_{MN},B_{MN},\Phi)$ in $10$-dimensional space-time,
supplemented by odd and even forms $C^{2k+1}$ and $C^{2k}$ respectively, together with self-duality constraints. Since we shall only discuss these theories tangentially in this paper, we omit the description of their actions.

\subsection{Supersymmetry}

We now say a few words about supersymmetry, which underlies much of what is discussed here. It is one of the main novel considerations beyond the equations suggested earlier by general relativity and gauge theories. Again, it is a deep concept, and we can only touch on a few of its mathematical consequences which motivate the equations we consider later.

\medskip
The key consequence of interest to us is that supersymmetry requires a spinor on space-time, which is covariantly constant with respect to a suitable connection which may have torsion.

\smallskip
This arises roughly as follows. Since we are dealing with supergravity theories, the bosonic fields always include a metric $G_{MN}$, and by supersymmetry, its partner $\chi_M{}$, which is a spinor-valued one-form $\chi_M$ called the gravitino field. Supersymmetry transformations act on the gravitino field as follows
\bea
\delta\chi_M=D_M \xi
\eea
where the infinitesimal generator $\xi$ is a spinor field and $D_M$ is a suitable connection on the spin bundle. Note that these transformations can be interpreted as the analogue of the infinitesimal variation of a metric $G_{MN}$ under diffeomorphisms generated by a vector field $V^M$, which is given by $\delta G_{MN}=\nabla_{\{M}V_{N\}}$.
Supersymmetry of a field configuration requires that, under a supersymmetry transformation, its gravitino field is unchanged. Thus we need a spinor $\xi$ satisfying $D_M\xi=0$.

Now we examine the possibilities for connections on the spin bundle. For our purposes, a spinor $\psi$ is just a section of a spinor bundle, and a spinor bundle is a vector bundle carrying a representation of the Clifford algebra $\{\gamma^M\}$, $\gamma^M\gamma^N+\gamma^N\gamma^M=2G^{MN}$. 
The simplest connection on spinors is the spin connection 
\bea
\nabla_M\xi=\partial_M\xi+{1\over 2}\omega_{MJN}\gamma^J\gamma^N\xi,
\eea 
where $\omega_{MJN}$ is the Levi-Civita connection. But other connections $D_M$ are possible, and may actually be required by the desired symmetries of the final theory. They necessarily differ from the spin connection by a Clifford algebra valued one-form and, expanding in the basis generated by $\gamma$-matrices, they are given by
\bea
D_M\xi
=
\nabla_M\xi+\sum_p\sum_{N_1,\cdots,N_p}H_{MN_1\cdots N_p}^{(p)}\gamma^{[N_1}\cdots\gamma^{N_p]}\xi
\eea
Cases of particular importance for $M$ theory are when the field $H_{MN_1\cdots N_p}^{(p)}$ is a $(p+1)$-form. The case of a $3$-form is responsible for the torsion $H$ in the equations for the Type I, Type II, and Heterotic string, and the case of a $4$-form results in the field $F_4$  in $11$-dimensional supergravity, both discussed earlier. In summary, we do find that supersymmetry of a configuration requires the existence of a spinor field $\xi$, covariantly constant  
with respect to a connection $D_M$ involving a torsion field $H$.

\medskip
We observe that the existence of a covariantly constant spinor is well known in mathematics to be characteristic of reduced holonomy and special geometry (see e.g. Berger \cite{Be}, Lichnerowicz \cite{Li}, Joyce \cite{J} and others). What unified string theories have provided is supersymmetry as a motivation, as well as the necessity of considering other connections differing from the Levi-Civita connection by a torsion field. More on the consequences of supersymmetry can be found in e.g. \cite{GGHPR, GMW, GP, GMR, Betal} and references therein.

\medskip
For phenomenological reasons, it is desirable to compactify space-time to $M^{3,1}\times X$, where $M^{3,1}$ is Minkowski or a maximally symmetric $4$-dimensional space-time, and $X$ is an internal space of dimension either $6$ or $7$, depending on whether we compactify string theories or $11$-dimensional supergravity. It is also desirable
to preserve supersymmetry upon compactification. The above considerations descend to similar considerations on the internal space $X$. In particular the existence of covariantly constant spinor fields $\xi$ imposes additional structure on the internal space, such as e.g. a $G_2$ structure when $X$ is $7$-dimensional or, when $X$ is $6$-dimensional, a complex structure and a holomorphic top-form $\Omega$, constructed from bilinears in $\xi$ as
\bea
J_M{}^N=i\xi^\dagger \gamma_M\gamma^N \xi,
\qquad
\Omega_{MNP}=\xi^\dagger \gamma_M\gamma_N\gamma_P\xi.
\eea
Thus we can trace back to supersymmetry the origin of the appearance of complex geometry in the equations of string theories.

\section{Equations from the Heterotic String}
\setcounter{equation}{0}

We come now to a detailed study of the equations resulting from the heterotic string. Here the constraints discussed earlier for supersymmetric compactifications have been worked out independently by C. Hull \cite{Hull2} and A. Strominger \cite{S}, who proposed the following system of equations.

\medskip

Let $X$ be a $3$-fold equipped with a holomorphic non-vanishing $(3,0)$-form $\Omega$ and a holomorphic vector bundle $E\to X$ with $c_1(E)=0$. Then the 
Hull-Strominger system is the following system for a Hermitian metric $\omega$ on $X$, with curvature $Rm\in \Lambda^{1,1}\otimes End(T^{1,0}(X))$, and a Hermitian metric $H_{\bar\alpha\beta}$ on $E$, with curvature $F\in \Lambda^{1,1}\otimes End(E)$, 
\bea
&&
i\partial\bar\partial \omega-{\alpha'\over 4}
{\rm Tr}(Rm \wedge Rm-F\wedge F)=0,
\qquad \omega^2\wedge F=0\nonumber\\
&&
d(\|\Omega\|_\omega\omega^2)=0
\eea
The last equation was written originally in terms of the torsion of $\o$. The above reformulation is due to Li and Yau \cite{LY} and will play a very important role below.
The Hull-Strominger system is a generalization of the solution found initially by Candelas, Horowitz, Strominger, and Witten \cite{CHSW}, which is obtained by setting $E=T^{1,0}(X)$ (up to a direct sum of flat bundles) and $H_{\bar\alpha\beta}=\omega$. Then the first and third equation reduce to $i\partial\bar\partial\omega=0$, $d(\|\Omega\|_\omega\omega^2)=0$, which imply that $\omega$ is K\"ahler and Ricci-flat. The second equation also reduces to the Ricci flat condition for $H_{\bar\alpha\beta}$, and so the system is consistent and all equations are satisfied. The emergence of K\"ahler Ricci-flat metrics, now known as Calabi-Yau metrics following the solution by Yau \cite{Y} of the Calabi conjecture, is the key unexpected link with complex geometry which we had mentioned earlier. We note that the second equation is just the Hermitian-Einstein equation, which is the bundle analogue of the K\"ahler Ricci-flat equation solved by Donaldson \cite{D} and Uhlenbeck-Yau \cite{UY}.

\subsection{The conformally balanced condition}

The Hull-Strominger system is a coupled system for a metric $\o$ on $X$ and a metric $H_{\bar\alpha\beta}$ on a vector bundle $E$. As a first step, we can consider a simplified situation where the metric $H_{\bar\alpha\beta}$ is known, and concentrate on the remaining two equations for $\o$. One of them, namely the third equation, is a $(2,2)$ cohomological condition, while the other, namely the first equation, is a condition on the curvature of a suitable representative in the cohomology class. 
These are typical features of canonical metrics in K\"ahler geometry, which are representatives of a fixed $(1,1)$ cohomology class satisfying e.g. the condition of having constant scalar curvature (see e.g. \cite{D1}, and for a survey \cite{PS}). From this point of view, the Hull-Strominger system can be viewed as providing a new notion of canonical metric in non-K\"ahler geometry.

\smallskip
Conditions of the form of the third equation had been introduced before in the mathematics literature. 
A metric $\omega$ is said to be {balanced} in the sense of Michelsohn \cite{M} is $\omega^2$ is closed. In this terminology, the third condition just says that the metric $\|\Omega\|_\o^{1/2}\o$ is balanced. The balanced condition is quite natural, and even has the advantage over the K\"ahler condition of being invariant under algebro-geometric modifications \cite{AB}. We shall say that $\o$ itself is conformally balanced. However, while formally similar, a $(2,2)$-form cohomological condition is much less manageable than a K\"ahler condition. A big difference is the {$\partial\bar\partial$-lemma} of K\"ahler geometry. It says that a metric $\omega$ is in the same K\"ahler class as $\omega_0$ if and only if
\bea
\omega=\omega_0+i\partial\bar\partial\varphi
\eea
where the potential $\varphi$ is unique up to a harmless additive constant. Any condition on the volume of $\o$ or a related curvature condition can then be easily rewritten as a partial differential equation in the potential $\varphi$. An example is the condition of vanishing Ricci curvature, which can be rewritten as a complex Monge-Amp\`ere equation on $\varphi$, and which has provided a particularly efficient way of finding and studying Calabi-Yau metrics.

\smallskip
By contrast, the absence of a $\p\bar\p$-lemma in non-K\"ahler geometry has been a major impediment in finding any efficient parametrization of balanced Hermitian metrics $\o$ with $\o^2$ in a given $(2,2)$-cohomology class, and hence in solving the Hull-Strominger system.

\subsection{Anomaly flows as substitutes for the $\partial\bar\partial$-lemma}

\medskip
We describe now joint works with Teng Fei, Sebastien Picard, and Xiangwen Zhang, based on the central idea that the absence of a $\partial\bar\partial$-lemma for balanced metrics can be bypassed by using a geometric flow. This idea was first applied to the case of the Hull-Strominger system in \cite{PPZ1,PPZ2,PPZ3, PPZ7}, and applied since more broadly to a range of other systems, including systems arising from the Type II A and Type II B strings \cite{FP, FPPZ} . 

\medskip
To be specific, consider the set-up for the Hull-Strominger system, that is, a compact $3$-fold $X$ equipped with a nowhere vanishing $(3,0)$ form $\Omega$
and a holomorphic vector bundle $E\to X$. We assume that $c_1(E)=0$ for notational simplicity. Let $\o_0$ be a conformally balanced Hermitian metric on $X$, in the sense that $d(\|\Omega\|_{\o_0}\o_0^2)=0$.
Then the Anomaly flow is defined to be the following flow for a pair $(\o(t),H_{\bar\alpha\beta}(t))$, where $\o(t)$ is a Hermitian metric on $X$ and $H_{\bar\alpha\beta}(t)$ a Hermitian metric on $E$,
\bea
\label{anomaly-HS}
&&
\partial_t(\|\Omega\|_\omega\omega^2)
=
i\partial\bar\partial\omega-{\alpha'\over 4}
{\rm Tr}\,(Rm\wedge Rm-F\wedge F)\nonumber\\
&&
H^{-1}\p_tH={\o^2\wedge F\over \o^3}.
\eea
Clearly the stationary points satisfy two of the three equations in the Hull-Strominger system. But the key point is that the second equation, namely the conformally balanced condition for $\o$, is automatically satisfied without having to appeal to a $\p\bar\p$-lemma. This is because the right hand side of the first equation in (\ref{anomaly-HS}) is closed, by standard Bott-Chern theory, and thus the closedness of $\|\Omega\|_\o\o^2$ is preserved along the flow.

\medskip
The flow (\ref{anomaly-HS}) admits many natural variants and generalizations, depending on the dimension, and on whether the well-known flow for the metric $H_{\bar\alpha\beta}$ can be decoupled. For our purposes, on complex manifolds $X$ of dimension $m\geq 2$, the variants which we consider fit into the following general flows for the sole metric $\o$,
\bea
\label{anomaly-Phi}
\p_t(\|\Omega\|_\o\o^{m-1})
=i\p\bar\p \o^{m-2}-\alpha'\Phi
\eea
where $\Phi=\Phi(\o, Rm(\o), t)$ is a closed $(m-1,m-1)$-form, and the initial data $\o(0)$ is conformally balanced in the sense that $d(\|\Omega\|_{\o_0}\o_0^{m-1})=0$.

\subsection{Formulation as a flow of metrics instead of $(2,2)$ forms}

While the Anomaly flow arises naturally as a flow of $(m-1,m-1)$-forms, its analysis requires an explicit and equivalent reformulation as a flow of $(1,1)$-forms. This was worked out in \cite{PPZ2, PPZ7}, using the explicit formulas derived there for the Hodge $\star$ operator. The answer is the following, where we have taken $\Phi=0$ for notational simplicity, the general case being obtainable in exactly the same way:

\begin{theorem} Consider the Anomaly flow with a {\it conformally balanced initial metric} on a complex manifold $X$ of dimension $m$. Then the flow can also be expressed as
\bea
\label{1-1flow}
\p_tg_{\bar kj}&=&{1\over (m-1)\|\Omega\|_\o}\big\{-\tilde R_{\bar kj}
-{1\over 2}T_{\bar kpq}\bar T_j{}^{pq}
+T_{\bar kjs}\bar \tau^s+\tau^{\bar r}\bar T_{j\bar k\bar r}+\tau_j\bar \tau_{\bar k}
\nonumber\\
&&
+{1\over 2(m-2)}(|T|^2-2|\tau|^2)g_{\bar kj}\big\}
\eea
for all $m\geq 3$. Here $\o=ig_{\bar kj}dz^j\wedge d\bar z^k$,
$\tilde R_{\bar kj}=g^{p\bar q}R_{\bar qp\bar kj}$ is the Chern-Ricci tensor, $T=i\partial\omega={1\over 2}T_{\bar kjm}dz^m\wedge dz^j\wedge d\bar z^k$ is the torsion tensor, and $\tau_\ell=(\Lambda T)_\ell=g^{j\bar k}T_{\bar kj\ell}$.
If $|\alpha' Rm(\omega(0))|$ is sufficiently small, the flow will exist at least for a short-time.
\end{theorem}

We note that, explicitly, 
$$
\tilde R_{\bar kj}=
-g^{p\bar q}g_{\bar km}\partial_{\bar q}( g^{m\bar \ell}\partial_pg_{\bar \ell j})
=-\Delta g_{\bar kj}+\cdots$$
which shows that the flow is parabolic for $\alpha'=0$. This parabolicity continues to hold for $\alpha'$ small enough. Note that the other Ricci tensor
$R_{\bar kj}=-\partial_j\partial_{\bar k}\log g=-g^{p\bar q}\partial_j\partial_{\bar k}g_{\bar qp}$ would not have been a good substitute for $\tilde R_{\bar kj}$ since it would not have implied parabolicity. When $m=3$, the expression between brackets on the right hand side of (\ref{1-1flow}) reduces to the simpler expression
$-\tilde R_{\bar kj}+g^{s\bar r}g^{p\bar q}T_{\bar qsj}\bar T_{p\bar r\bar k}$, which was the form of the flow originally derived in \cite{PPZ2}. That the two expressions actually coincide is a consequence of an identity for the torsion tensor found in \cite{FP}.

\smallskip
The original formulation of the Anomaly flow as well as the above formula for $\p_tg_{\bar kj}$ require the existence of the nowhere vanishing holomorphic form $\Omega$. Recently, in the case $\Phi=0$, another formulation was found where $\Omega$ appeared only in the initial data, and hence the flow can be generalized to manifolds where such a form $\Omega$ may not exist \cite{FP}:

\begin{theorem}
Let $\eta$ be the Hermitian metric defined by $\eta=\|\Omega\|_\o\o$. If $\o$ evolves by the Anomaly flow (\ref{anomaly-Phi}) with $\Phi=0$
and $\o(0)$ is conformally balanced in the sense that $d(\|\Omega\|_{\o(0)}\o^{m-1})=0$, then the corresponding metrics $\eta$ evolve according to the following flow
\bea
\label{anomaly-eta}
i^{-1}\p_t\eta=-{1\over m-1}(\tilde R_{\bar kj}(\eta)+{1\over 2}T_{\bar kpq}(\eta)\bar T_j{}^{pq}(\eta)).
\eea
and its initial data satisfies $d(\|\Omega\|_{\eta(0)}^2\eta^{m-1})=0$. In particular, rescaling time by $t\to (m-1)^{-1}t$, the Anomaly flow admits a natural extension to all dimensions.
\end{theorem}

This formulation allows a generalization of the Anomaly flow to arbitrary initial data, and not just conformally balanced ones. These generalizations fit into the family of generalizations of the Ricci flow introduced in \cite{ST}. The flow (\ref{anomaly-eta})  actually coincides with the flow identified in \cite{U} as a flow preserving the Griffiths positivity and the dual Nakano-positivity of the tangent bundle. Such properties hold for the K\"ahler-Ricci flow \cite{B, M}, and are known to have many important consequences, see e.g. \cite{Ha, PSSW}.

\subsection{The Anomaly flow and the Fu-Yau solution of the Hull-Strominger system}

The Anomaly flow method appears to have a lot of potential. Here we shall discuss how it can be applied to give new proofs of two fundamental results in complex geometry, namely the Fu-Yau solution of the Hull-Strominger system, and Yau's solution of the Calabi conjecture. We begin with the Fu-Yau solution.

\smallskip

Let $(Y,\hat\omega)$ be a Calabi-Yau surface, equipped with a nowhere vanishing holomorphic form $\Omega_Y$. Let $\omega_1,\omega_2\in H^2(Y,{\bf Z})$ satisfy $\omega_1\wedge\hat\omega=\omega_2\wedge\hat\omega=0$.

\medskip
From this data, building on earlier ideas of Calabi and Eckmann \cite{CE}, Goldstein and Prokushkin \cite{GP} constructed
a toric fibration $\pi:X\to Y$, equipped with a $(1,0)$-form $\theta$ on $X$ satisfying
$\partial\theta=0$, $\bar\partial\theta=\pi^*(\omega_1+i\omega_2)$. Furthermore,
the form 
\bea
\Omega=\sqrt 3\,\Omega_M\wedge\theta
\eea
is a holomorphic nowhere vanishing $(3,0)$-form on $X$, and for any scalar function $u$ on $Y$, the $(1,1)$-form
\bea
\omega_u=\pi^*(e^u\hat\omega)+i\theta\wedge\bar\theta
\eea
is a conformally balanced metric on $X$. In a breakthrough on the Hull-Strominger system, Fu and Yau \cite{FY1, FY2} found the first non-perturbative, non-K\"ahler solution by showing how the system would descend on these fibrations to a single scalar Monge-Amp\`ere equation. They succeeded in solving this equation, even though it involved new difficult gradient terms. Using the Anomaly flow, we find \cite{PPZ3}:

\begin{theorem}
\label{anomaly-FY}
Consider the Anomaly flow on the fibration $X\to Y$ constructed above,
with initial data $\omega(0)=\pi^*(M\hat\omega)+i\theta\bar\theta$, where $M$ is a positive constant. Then for a suitable bundle $\pi^*(E_Y)$ on $X$ and a suitable initial Hermitian metric $H_{\bar\alpha\beta}(0)$ (see below), the metrics
$\omega(t)$ are of the form $\pi^*(e^u\hat\omega)+i\theta\wedge \bar\theta$. Assuming an integrability condition on the data (which is necessary), there exists $M_0>0$, so that for all $M\geq M_0$, the flow exists for all time, and converges to a metric $\omega_\infty$ with $(\omega_\infty, H_{\bar\alpha\beta})$ satisfying the Hull-Strominger system.
\end{theorem}

We note that, in general, a given elliptic equation admits an infinite number of parabolic flows with it as stationary point. But not all of these may have the right structural properties to satisfy long-time existence and to converge. The Anomaly flow was motivated by a completely different, geometric, consideration, namely to preserve the conformally balanced condition. The above theorem on the present case of toric fibrations suggests that it is also well-behaved from the analytic viewpoint.

\medskip
We describe some of the main ideas in the proof of Theorem \ref{anomaly-FY}. As in \cite{FY1, FY2}, if we choose the bundle $E$ to be the pull-back $\pi^*(E_Y)$ of a bundle $E_Y$ on the base $Y$, then the system can be shown to descend to a system on the base $Y$,
for an unknown metric of the form $\hat\o_u=e^u\hat\omega$. We can take $H$ to be the pull-back of a metric on $E_Y$ which is Hermitian-Einstein with respect to $\hat\o$, and hence with respect to $\hat\o_u$ for any $u$. Then the Anomaly flow descends to the following flow on the base $Y$,
\bea
\label{anomaly-baseY}
\partial_t\hat\omega_u
=
-{1\over 2\|\Omega\|_{\hat\omega_u}}({R(\hat\omega_u)\over 2}-|T(\hat\omega_u)|^2
-
{\alpha'\over 4}\sigma_2(iRic_{\hat\omega(u)})+2\alpha'
{i\partial\bar\partial (\|\Omega\|_{\hat\omega_u} \rho)\over \hat\omega_u^2}-2{\mu\over\hat\omega_u^2})\hat\omega_u
\nonumber
\eea
where $\rho$ and $\mu$ are fixed and given forms which depend on the geometric data $(Y,\Omega_Y,\omega_1,
\omega_2)$.

\medskip
To start, we assume that the initial data satisfies $|\alpha' Ric_\omega|<<1$, so that the diffusion operator
$$
\Delta_F=F^{p\bar q}\nabla_p\nabla_{\bar q},
\qquad
F^{p\bar q}=g^{p\bar q}+\alpha'\|\Omega\|_\omega^3\tilde\rho^{p\bar q}
-
{\alpha'\over 2}(R g^{p\bar q}-R^{p\bar q})
$$
is positive definite. The key and difficult step is to prove that this condition is preserved along the flow. This is done through the following several careful estimates in terms of the scale $M$:

\medskip

$\bullet$ {\it Uniform equivalence of the metrics $\omega(t)$}:
this is equivalent to a uniform estimate for the conformal factor $u$ in $\hat\omega_u=e^u\hat\omega$, and is established by Moser iteration, exploiting the fact that the quantity $\int_X \|\Omega\|_\omega\omega^2$ is conserved along the flow. For our purposes, we shall need the following precise version in terms of $M$.
Assume that the flow exists for $t\in [0,T)$ and starts with $\hat\omega(0)=M\hat\omega$. Then there exists $M_0$ so that, for $M\geq M_0$, we have
\bea
{\rm sup}_{X\times [0,T)}e^u\leq C_1M,
\qquad
{\rm sup}_{X\times [0,T)}e^{-u}\leq {C_2\over M}
\eea
where $C_1,C_2$ depend only on $(X,\hat\omega)$, $\mu$, $\rho$, and $\alpha'$.

\medskip
$\bullet$ {\it Estimates for the torsion}:
there exists $M_0$ with the following property. If the flow is started with $\omega(0)=M\hat\omega$ and $M\geq M_0$, and if
\bea
|\alpha' Ric_{\hat\omega}|\leq 10^{-6}
\eea
along the flow, then there exists a constant $C_3$ depending only on $Y,\hat\omega,\mu, \rho,\alpha'$ so that
\bea
|T|^2\leq {C_3\over M^{4/3}}<<1.
\eea

\medskip

$\bullet$ {\it Estimates for the curvature}:
start the flow with $\omega(0)=M\hat\omega$. There exists $M_0>1$ such that, for every $M\geq M_0$, if
{$$
\|\Omega\|^2\leq {C_2^2\over M^2},
\qquad 
|T|^2\leq {C_3\over M^{4/3}}
$$}
along the flow, then
{$$
|\alpha'\,Ric_{\hat\omega}|\leq {1\over M^{1/2}}
$$}

\medskip

$\bullet$ {\it Higher order estimates}:
under similar conditions, we can establish estimates for the higher order derivatives of the torsion and the curvature. An interesting new technical point is the usefulness of test functions of the form
{$$
G=(|\alpha' Ric_{\hat\omega}|+\tau_1)
|\nabla Ric_{\hat\omega}|^2
+
(|T|^2+\tau_2)|\nabla T|^2
$$}

\smallskip

$\bullet$ {\it Closing the loop of estimates}:
if we start with an initial data satisfying $|\alpha' Ric_{\hat\omega}|\leq 10^{-6}$, the estimates for the torsion imply that $|T|^2\leq C_3M^{-4/3}$ and hence $|\alpha' Ric_{\hat\omega}|\leq M^{-1/2}$.  But this implies that for $M>>1$, the flow does not leave the region $|\alpha' Ric_{\hat\omega}|\leq 10^{-6}$, and hence all these estimates hold for all time. This implies the existence for all time of the flow, and an additional argument establishes its convergence in $C^\infty$.

\medskip
We observe that, as anticipated, the techniques discovered in the solution of the Anomaly flow have proven to be useful for other partial differential equations as well \cite{PPZ4, PPZ5,PPZ6}. Earlier techniques can be found in \cite{Guanbo, GRW}.

\subsection{The Anomaly flow and the Calabi conjecture}

Next, we apply the Anomaly flow to obtain a new proof, with new estimates, of Yau's theorem on the existence of Ricci-flat K\"ahler metrics. For this it suffices to consider the special case with $\alpha'=0$ in the flow (\ref{anomaly-Phi}). It is not difficult to show that conformally balanced metrics $\o$ which also satisfy the stationary condition of this flow, namely $i\p\bar\p\o=0$, must be K\"ahler and Ricci-flat metrics \cite{CHSW, MT, FT, PPZ7}. Thus it suffices to produce an initial data for which the flow (\ref{anomaly-Phi}) with $\alpha'=0$ converges. In \cite{PPZ7}, we prove:

\begin{theorem}
Assume that the initial data $\omega(0)$ satisfies
$\| \Omega \|_{\omega(0)} \omega(0)^{n-1} = \hat{\chi}^{n-1}$,
where $\hat{\chi}$ is a K\"ahler metric. Then the flow exists for all time $t>0$, and as $t \rightarrow \infty$, the solution $\omega(t)$ converges smoothly to a K\"ahler, Ricci-flat, metric $\omega_\infty$. If we define the metric $\chi_\infty$ by
\bea
\omega_\infty = \| \Omega \|_{\chi_\infty}^{-2/(n-2)} \chi_\infty,
\eea
then $\chi_\infty$ is the unique K\"ahler Ricci-flat metric in the cohomology class $[\hat{\chi}]$, and $\| \Omega \|_{\chi_\infty}$ is an explicit constant.
\end{theorem}

\medskip
In this case, the Anomaly flow is actually conformally equivalent to a flow of metrics $t\to \chi(t)$ in the K\"ahler class of $\hat\chi$. Indeed, given a solution $\o(t)$ of the Anomaly flow, we can define the Hermitian metric $\chi(t)$ by
\bea
\label{anomaly-chi}
\|\Omega\|_{\o(t)}\o^{m-1}(t)=\chi^{m-1}(t).
\eea
Then the flow (\ref{anomaly-Phi}) with $\alpha'=0$ can be rewritten as
\bea
\p_t\chi^{m-1}(t)=i\p\bar\p(\|\Omega\|_{\chi(t)}^{-2}\chi^{m-2}(t))
\eea
and hence, upon carrying out the differentiations and assuming that $\chi(t)$ is K\"ahler,
\bea
(m-1)\p_t\chi\wedge\chi^{m-2}=(i\p\bar\p\|\Omega\|_{\chi(t)}^{-2})\wedge\chi^{m-2}.
\eea
This last equation is satisfied if $\chi(t)$ flows according to
\bea
\label{anomaly-chi1}
(m-1)\p_t\chi=i\p\bar\p\|\Omega\|_{\chi(t)}^{-2}.
\eea
But this flow preserves the K\"ahler condition, and hence all the steps of the previous derivation are justified. In particular, by the uniqueness of solutions of parabolic equations, the solution $\o(t)$ of (\ref{anomaly-Phi}) is indeed given by (\ref{anomaly-chi}) with $\chi(t)$ K\"ahler.

\smallskip
The flow (\ref{anomaly-chi1}) can be written explicitly in terms of potentials. Setting
$\chi(t) = \hat{\chi} + i \ddb \varphi(t) > 0$, we find
\bea
\p_t \varphi = e^{-f} {(\hat\chi+i\partial\bar\partial\varphi)^n
\over\hat\chi^n}, \ \ \varphi(z,0)=0
\eea
where the scalar function $f\in C^\infty(X,{\bf R})$ is defined by
$(n-1)e^{-f} = \| \Omega \|^{-2}_{\hat{\chi}}$.
We note that this flow is of Monge-Amp\`ere type, but without the log as in the solution by K\"ahler-Ricci flow by Cao \cite{Cao}, and without the inverse power of the determinant, as in the recent equation proposed by Cao and Keller \cite{CaK} and Collins, Hisamoto, and Takahashi \cite{CHT}. Because of this, we need a new way of obtaining $C^2$ estimates. It turns out that the test function 
$$
G(z,t)=\log {\rm Tr}\, h-A(\varphi-{1\over[\hat\chi^n]}\int_X\varphi\hat\chi^n)+B[{(\hat\chi+i\partial\bar\partial\varphi)^n\over\hat\chi^n}]^2
$$
can do the job. We observe that it differs from the standard test function $\hat G(z,t)
=\log {\rm Tr}\, h-A\varphi$ used in the study of Monge-Amp\`ere equations (see e.g. \cite{PSS} and references therein) by terms involving the square of the Monge-Amp\`ere determinant. Thus the study of the Anomaly flow leads again to new tools that should be useful in the development of the theory of non-linear partial differential equations.

\subsection{Further remarks}

While the Anomaly flow has been very successful in the several cases studied so far, its general theory remains to be more fully developed. From the sole fact that it can be viewed as an extension of the K\"ahler-Ricci flow to the non-K\"ahler case, the general theory can be expected to be complicated. We restrict ourselves to a few remarks.

\medskip
It may be instructive to draw a closer comparison of the Anomaly flow with the K\"ahler-Ricci flow.
In the K\"ahler-Ricci flow, the $(1,1)$-cohomology class is determined
\bea
[\omega(t)]=[\omega(0)]-tc_1(X)
\eea 
In the Anomaly flow, we have something similar
\bea
[\|\Omega\|_{\omega(t)}\omega(t)^2]
=
[[\|\Omega\|_{\omega(0)}\omega(0)^2]-t\alpha'(c_2(X)-c_2(E)).
\eea
However, the $(2,2)$ cohomology class $[\|\Omega\|_{\omega}\omega^2]$ provides much less information than the $(1,1)$-cohomology class $[\omega]$. For example, the volume is an invariant of an $(1,1)$ cohomology class, but not of an $(2,2)$-cohomology class.
As an indirect consequence, the maximum time of the Anomaly flow is not determined by cohomology alone and depends on the initial data. In this delicate dependence on the initial data, the Anomaly flow is closer to the Ricci flow than the K\"ahler-Ricci flow.

\medskip
The dependence of the long-time behavior of Anomaly flow on the initial data can be seen explicitly in several examples worked out by Fei, Huang, and Picard \cite{FHP1, FP}, including on hyperk\"ahler fibrations over Riemann surfaces. It is shown there that the flow can both exist for all time or terminate in finite time. In the first case, the flow, suitably normalized, can collapse the fibers. But even in this very specific geometric setting the behavior of the flow has not been worked out for general data. In particular, it has not been worked out for data close to the stationary points found earlier in \cite{FHP2}, which are particularly interesting as an infinite family of topologically distinct solutions to the Hull-Strominger system. That such a family of Calabi-Yau metrics is conjectured, but not yet known, to exist is an illustration of the comparative flexibility of the Hull-Strominger system.

\medskip
A priori, one of the reasons the Anomaly flow may terminate is if $\|\Omega\|_\o$ goes to $0$ or $\infty$. But at least when $\alpha'=0$, it has now been shown that $\|\Omega\|_\o$ always remains bounded along the flow \cite{FP}. An essential tool is the monotonicity of dilaton functionals, which had also been considered in \cite{GF1, GF2, GFRST}.

\smallskip
Finally, it would be helpful to work out many more examples. In this context, we note that many solutions of the Hull-Strominger system have been found in many geometric settings, which should be instructive to investigate, see e.g. \cite{FIUV, FIUV2, FGV, GF1, GF2, GFRST} and references therein.

\section{Equations from the Type II A and Type II B strings}
\setcounter{equation}{0}

Next, we consider supersymmetric compactifications of the Type II B string and of the Type II A string with brane sources, as formulated by Tseng and Yau \cite{TY1}, building on earlier formulations of Grana-Minasian-Petrini-Tomasiello \cite{Gr}, Tomasiello \cite{T}, and others. In particular, (subspaces of) linearized solutions have been identified by Tseng and Yau with Bott-Chern and Aeppli cohomologies in the case of Type II B, and with their own symplectic cohomology in the case of Type II A, as well as interpolating cohomologies \cite{TY2,TY3}. Related boundary value problems have been recently studied by Tseng and Wang \cite{TW}.

\medskip
The main feature of these equations of concern to us is that they all involve cohomological conditions that are not known to be enforceable by any $\p\bar\p$-lemma.
Here we shall just point out that the same idea of Anomaly flows can be applied. A more detailed study of the resulting flows will appear elsewhere
\cite{FPPZ}. 
\medskip

\subsection{Type II B strings}

\smallskip

Let $X$ be compact $3$-dimensional complex manifold, equipped with a nowhere vanishing holomorphic $3$-form $\Omega$. Let $\rho_B$ be the Poincare dual of a linear combination of holomorphic $2$-cycles. We look for a Hermitian metric satisfying the following system
{$$
d\omega^2=0, \qquad i\partial\bar\partial (\|\Omega\|_\omega^{-2}\omega)=\rho_B$$}
where $\|\Omega\|_\omega$ is defined by $i\Omega\wedge\bar\Omega=\|\Omega\|_\omega^2\omega^3$. If we set $\eta=\|\Omega\|_\omega^{-2}\omega$, this system can be recast in a form similar to the Hull-Strominger system,
{$$
d(\|\Omega\|_\eta\eta^2)=0,\qquad i\partial\bar\partial \eta=\rho_B$$}
This system can be approached by the following Anomaly type flow,
\bea
\partial_t(\|\Omega\|_\eta\eta^2)=i\partial\bar\partial\eta-\rho_B,
\quad
d(\|\Omega\|_{\eta_0}\eta_0^2)=0.
\eea
Because the right hand side is closed, the closedness of the initial condition is preserved, and the system is solved if the flow converges.

\subsection{Type II A strings}

\smallskip

Let $X$ be this time a real $6$-dimensional symplectic manifold, in the sense that it admits a closed, non-degenerate $2$-form $\omega$ (but there may be no compatible complex structure, so it may not be a K\"ahler form). Then the equations are now for a complex $3$-form $\Omega$ with ${\rm Im} \,\Omega=\star{\rm Re}\,\Omega$, and
{
$$d({\rm Re}\,\Omega)=0,
\qquad
dd^\Lambda(\star\|\Omega\|^2{\rm Re}\,\Omega)=\rho_A$$}
where $\rho_A$ is the Poincare dual of a linear combination of special Lagrangians. Here $d^\Lambda=d\Lambda-\Lambda d$ is the symplectic adjoint. 
Thus it is natural to introduce another Anomaly-type flow
\bea
\partial_t({\rm Re}\,\Omega)=dd^\Lambda(\star \|\Omega\|^2{\rm Re}\,\Omega)-\rho_A,
\quad
d({\rm Re}\,\Omega_0)=0.
\eea
whose stationary points would again solve the desired system.

\section{Equations from $11$-dimensional Supergravity}
\setcounter{equation}{0}

It does not appear that the mathematical study of the equations from $11$-dimensional supergravity is as extensive as in the heterotic case. Nevertheless we shall describe joint works with Teng Fei and Bin Guo on some exact solutions and their moduli, criteria for the construction of supersymmetric compactifications, and the formulation of some geometric flows which may provide an approach to more general solutions \cite{FGP1, FGP2}.

\medskip

Recall that the fields of $11$-dimensional supergravity are an $11$-dimensional Lorentz metric $G_{ij}$ and a $4$-form $F=dA$, and the action is given by (\ref{11d-action}).
The resulting field equations are
$$
d\star F={1\over 2}F\wedge F,
\qquad
R_{ij}={1\over 2}(F^2)_{ij}-{1\over 6}|F|^2 G_{ij}
$$
where the symmetric 2-tensor $F^2$ is defined by
$$
(F^2)_{ij}={1\over 6}F_{iklm}F_j{}^{klm}.$$
The supersymmetric solutions are the solutions which admit a spinor $\xi$ satisfying
$$
D_m\xi:=\nabla_m\xi-{1\over 288}F_{abcd}(\Gamma^{abcd}{}_m+8\Gamma^{abc}\delta^d{}_m)\xi=0$$
i.e. spinors which are covariantly constant with respect to the connection $D_m$, obtained by twisting the Levi-Civita connection with the flux $F$.

\subsection{Early solutions}

\smallskip
Some early solutions were found with the Ansatz $M^{11}=M^4\times M^7$, where $M^4$ is a Lorentz 4-manifold and $M^7$ a Riemannian manifold with metrics $g_4$ and $g_7$ respectively. Setting $F=c Vol_4$ where $Vol_4$ is the volume form on $M^4$ reduces the field equations to
$$
(Ric_4)_{ij}=-{c^2\over 3}(g_4)_{ij},\qquad
(Ric_7)_{ij}={c^2\over 6}(g_7)_{ij}$$
i.e., $M^4$ and $M^7$ are Einstein manifolds with negative and positive scalar curvatures respectively. These are the Freund-Rubin solutions, which include $AdS^4\times S^7$ \cite{FR}.

\smallskip
More sophisticated solutions can be found with other ansatz for $F$, e.g. $F=cVol_4+\psi$ for suitable $\psi$, leading to nearly $G_2$ manifolds, as well as many others see e.g. \cite{E, PvN, PW, DH, Du}.

\medskip

For us the special solution of particular interest, following Duff-Stelle \cite{DS}, is obtained by setting $M^{11}=M^3\times M^8$, 
$$
g_{11}=e^{2A}g_3+g_8,
\qquad
F=Vol_3\wedge df$$
where $g_3$ is a Lorentz metric on $M^3$, $g_8$ is a Riemannian metric on $M^8$, and $(A,f)$ are smooth functions on $M_8$. The now well-known solution of Duff-Stelle is then obtained by assuming the flatness of $g_3$, the conformal flatness of $g_8$, the radial dependence of $A,f$, and supersymmetry.

\subsection{Other multimembrane solutions}

We now discuss a way of finding more systematically solutions to $11$-dimensional supergravity \cite{FGP1}. To begin with, we consider solutions given by warped products $M^{11}=M^3\times M^8$, $g_{11}=e^{2A}g_3+g_8,
F=Vol_3\wedge df$ as in the original work of Duff and Stelle. The first result is a complete characterization of such data giving rise to a supersymmetric solution:

\begin{theorem}
The data $(g_3,g_8,A,f)$ is a supersymmetric solution to 11-dimensional supergravity equation if and only if

{\rm (a)} $g_3$ is flat;

{\rm (b)} $\bar g_8 := e^Ag_8$ is a Ricci-flat metric admitting covariantly constant spinors with respect to the Levi-Civita connection;

{\rm (c)} $e^{-3A}$ is a harmonic function on $(M_8,g_8)$ with respect to the metric $\bar g_8$; 

{\rm (d)} $df = \pm d(e^{3A})$.
\end{theorem}

Applying the classic results of S.Y. Cheng, P. Li, and S.T. Yau \cite{CL, LpY} on lower bounds for Green's functions as well as new constructions of complete K\"ahler-Ricci flat metrics on ${\bf C}^4$ by Szekelyhidi \cite{Sze}, Conlon-Rochon \cite{CR}, and Li \cite{Li}, we obtain in this manner the complete supersymmetric solution for this ansatz, which includes many solutions not available before in the physics literature.

\medskip

Next, we find indications that $11$-dimensional supergravity may have some integrable structure. Indeed, we find that the solution found by Duff-Stelle can be imbedded into a five-parameter family of solutions, the only one which is supersymmetric is the Duff-Stelle solution itself. More precisely \cite{FGP1},

\begin{theorem}
There is a 5-parameter family of solutions $(g_3,g_8,A.f)$ of solutions of the field equations of $11$-dimensional supergravity. In fact, $g_3$ is Einstein and one of the parameters is its scalar curvature, a second parameter is a global scale parameter for $g_8$, and the three remaining parameters are constants of integration for the following 3rd order ODE, from which the rest of the Ansatz can be determined,
\bea
{d^3v\over dt^3}+7{d^2v\over dt^2}v+14 ({dv\over dt}^2)+2{dv\over dt}(17v^2-60)
+12(v^2-4)(v^2-6)=0.
\eea
\end{theorem}

\bigskip
Actually many special solutions of this equation can be written down explicitly, and it may deserve further investigation from the pure ODE viewpoint.

\subsection{Parabolic reductions of $11$-dimensional supergravity}

More general solutions of $11$-dimensional supergravity will ultimately have to be found by solving partial differential equations. These equations will be hyperbolic, because of the Lorentz signature of $M^{11}$, but as a start, we can try and identify the subcases where the Lorentz components are known, and deal only with elliptic equations. Thus
we consider space-times and field configurations of the form $M^{11}=M^{1,p}\times M^{10-p}$,
\bea
G=e^{2A}g_{1,p}+g,
\quad F=dVol_g\wedge\beta+\Psi
\eea
where $\beta$ and $\Psi$ are now respectively a $1$-form and a $4$-form on $M^{10-p}$, both closed \cite{FGP2}:

\begin{theorem}
There exist general parabolic flows of the configuration $(g(t), A(t),\beta(t),\Psi(t))$ (which can be written down explicitly) with the following properties:

{\rm (a)} The forms $\beta$ and $\Psi$ remain closed along the flow, and the above ansatz is preserved;

{\rm (b)} The corresponding configuration $(G,F)$ on $M^{11}$ evolves in time by the flow
$$
\partial_t G_{MN}=-2R_{MN}+F_{MN}^2-{1\over 3}|F|^2G_{MN},
\quad
\partial_tF=-\Box F-{1\over 2}d\star(F\wedge F)$$
whose stationary points are (assuming a cohomological condition) solutions of $11$-dimensional supergravity;

{\rm (c)} If $T<\infty$ is the maximum time of existence of the flow, then
$$
{\rm limsup}_{t\to T^-}{\rm sup}_{M^{10-p}}
(|Rm|+|A|+|\beta|+|\Psi|)=\infty$$

\end{theorem}

The constraint of supersymmetry has not yet been enforced on these reductions, and this is an important avenue for further investigations. We note in this context the wealth of results for supersymmetric solutions of supergravity theories in dimensions 5 and 6 \cite{GGHPR, GMW, GP, GMR, Betal}.
In a different vein, interesting related flows of balanced metrics or $G_2$ structures have been proposed in \cite{BV, BryXu, Ka, LoWe, FFR}. Yet another line of investigation is that of static solutions with symmetry, pioneered in \cite{We} for Einstein's equations in 4 dimensional, and more recently in \cite{AKK} for $5$-dimensional minimal supergravity.

\bigskip
{\bf Acknowledgements} The author would like to thank the organizers of the ICCM 2018 in Taipei, Taiwan, for their invitation to speak there. Part of this material was also presented at the Johns Hopkins conference in honor of B. Shiffman,
at a Colloquium at the University of Connecticut, Storrs, and at the conference in honor of S.T. Yau at Harvard University in May 2019.
The author would like to thank the Galileo Galilei Institute for Theoretical Physics and INFN for hospitality and partial support during the workshop ``String Theory from a worldsheet perspective" where part of this paper was written. He would also like to thank Professor Rodolfo Russo and Professor Li-Sheng Tseng for very stimulating conversations and correspondence.

\bigskip

email address: phong@math.columbia.edu

Department of Mathematics

Columbia University, New York, NY 10027 USA

\end{document}